\documentclass{article}
\catcode`\@=11
\newtheorem{theorem}{Theorem}
\newtheorem{lemma}{Lemma}
\newtheorem{proposition}{Proposition}

\newtheorem{example}{Example}

\newtheorem{corollary}{Corollary}
\def\demo{\noindent{\bf Proof .-}}
\def\section{\@startsection {section}{1}{\z@}{-3.5ex plus -1ex
minus-.2ex}{2.3ex plus .2ex}{\normalsize\bf}}

\def\bz{\hbox{\it Z\hskip -4pt Z}}

\newcommand{\het}{H_{\rm et}}
\newcommand{\hc}{H_{\rm c}}
\begin{document}
\begin{center}
{\Large\bf \textsc{On ideals generated by monomials and one binomial}}
\end{center}
\vskip.5truecm
\begin{center}
{Margherita Barile\footnote{Partially supported by the Italian Ministry of Education, University and Research.}\\ Dipartimento di Matematica, Universit\`{a} di Bari, Via E. Orabona 4,\\70125 Bari, Italy}\footnote{e-mail: barile@dm.uniba.it, Fax: 0039 080 596 3612}
\end{center}
\vskip1truecm
\noindent
{\bf Abstract} We determine, in a polynomial ring over a field, the arithmetical rank of certain ideals generated by a set of monomials and one binomial.  

\section*{Introduction} Let $R$ be Noetherian commutative ring with identity. We say that some elements  $p_1,\dots, p_s\in R$ generate an ideal $I$ of $R$  {\it up to radical} if $\sqrt{I}=\sqrt{(p_1,\dots, p_s)}$. The smallest $s$ with this property is called the {\it arithmetical rank} of $I$, denoted ara\,$I$. This number is of great interest in Algebraic Geometry: if  $R$ is a polynomial ring in $n$ indeterminates over an algebraically closed field $K$,  by virtue of Hilbert Basissatz, it is the minimum number of defining equations for the variety $V(I)$ in the affine space $K^n$.\newline It is well-known that  height\,$I\leq$\,ara\,$I$; if equality holds, $I$ is called a {\it set-theoretic complete intersection} (s.t.c.i.). Obviously, ara\,$I\leq\mu(I)$, where $\mu(I)$ is the minimum number of generators for $I$. This inequality is strict in general, which implies that ara\,$I$ elements generating $I$ up to radical need to be constructed by combining the minimal generators in a suitable way. No universal method is known so far, but some techniques have been developed for special cases. One example is the following lemma by Schmitt and Vogel:
\begin{lemma}\label{lemma1}{\rm [\cite{SV}, p.~249]} Let $P$ be a finite subset of elements of $R$. Let $P_0,\dots, P_r$ be subsets of $P$ such that
\begin{list}{}{}
\item[(i)] $\bigcup_{l=0}^rP_l=P$;
\item[(ii)] $P_0$ has exactly one element;
\item[(iii)] if $p$ and $p''$ are different elements of $P_l$ $(0<l\leq r)$ there is an integer $l'$ with $0\leq l'<l$ and an element $p'\in P_{l'}$ such that $pp''\in(p')$.
\end{list}
\noindent
We set $q_l=\sum_{p\in P_l}p^{e(p)}$, where $e(p)\geq1$ are arbitrary integers. We will write $(P)$ for the ideal of $R$ generated by the elements of $P$.  Then we get
$$\sqrt{(P)}=\sqrt{(q_0,\dots,q_r)}.$$

\end{lemma}
This result, together with the extensions and refinements contained in \cite{B1}, \cite{B2} and \cite{B3}, works for many monomial ideals. The propositions  in \cite{B0}, \cite{BMT}, \cite{Br}, \cite{BS}, \cite{El}, \cite{E},  and  \cite{RV} apply to important classes of binomial ideals, namely to determinantal and toric ideals. The aim of this paper is to  present new methods for the computation of the arithmetical rank of certain polynomial ideals which have a ``mixed" generating set, formed by {\it one} binomial and some monomials.  We shall  compare their arithmetical ranks with the lower bounds provided by local and \'etale cohomology. The case of the ideals defining minimal varieties, which are generated by {\it several} binomials and monomials, has been treated in \cite{B4}. 
\section{Generalizing a result by Schmitt and Vogel}
We prove the following generalization of Lemma \ref{lemma1}. The proof is an easy adaptation of the original one. We give it for the sake of completeness.
\begin{lemma}\label{lemma2}Let $P$ be a finite subset of elements of $R$. Let $P_0,\dots, P_r$ be subsets of $P$ such that
\begin{list}{}{}
\item[(i)] $\bigcup_{l=0}^rP_l=P$;
\item[(ii)] $P_0$ has exactly one element;
\item[(iii)] if $p$ and $p'$ are different elements of $P_l$ $(0<l\leq r)$ then $(pp')^m\in\left(\bigcup_{i=0}^{l-1}P_i\right)$ for some positive integer $m$. 
\end{list}
\noindent
We set $q_l=\sum_{p\in P_l}p^{e(p)}$, where $e(p)\geq1$ are arbitrary integers. Then we get
$$\sqrt{(P)}=\sqrt{(q_0,\dots,q_r)}.$$
\end{lemma}
\demo The inclusion $\supset$ is clear. For $\subset$ it suffices to prove that, for $0\leq l\leq r$, $P_l\subset \sqrt{(q_0,\dots, q_l)}$. We proceed by induction on $l$. For $l=0$ the claim is trivial, since by (ii) $P_0=\{q_0\}$. So suppose that $l\geq1$ and assume the claim true for $l-1$. Let $P_l=\{p_0,\dots, p_s\}$, and consider
$$p_0q_l=p_0^{e(p_0)+1}+\sum_{i=1}^sp_0p^{e(p_i)}.$$
\noindent
Then take the $N$-th power of both sides, where $N$ is a positive integer. One gets:
\begin{equation}\label{1}(p_0q_l)^N=p_0^{N(e(p_0)+1)}+C,\end{equation}
\noindent
where $C$ is a linear combination of products of the form 
$$p_0^{(e(p_0)+1)\alpha_0+\alpha_1+\cdots+\alpha_s}p_1^{e(p_1)\alpha_1}\cdots p_s^{e(p_s)\alpha_s},$$
\noindent
 with $\sum_{i=0}^s\alpha_i=N$. For $N$ large enough, from (iii) and induction it follows that 
\begin{equation}\label{1'}C\in\left(\bigcup_{i=0}^{l-1}P_i\right)\subset \sqrt{(q_0,\dots,q_{l-1})}.\end{equation}
Then (\ref{1}) and  (\ref{1'}) imply that
$$p_0^{N(e(p_0)+1)}=(p_0q_l)^N-C\in\sqrt{(q_0,\dots, q_{l-1}, q_l)},$$
\noindent 
whence
$$p_0\in \sqrt{(q_0,\dots,q_{l-1}, q_l)}.$$
\noindent
Since $p_0$ is an arbitrary element of $P_l$, this completes the induction step.
\par\medskip\noindent
\begin{example}\label{example1}{\rm Let $K$ be a field and consider the polynomial ideal $I\subset R=K[x_1,\dots, x_6]$ generated by the set
$$P=\{x_1x_2+x_3x_4, x_1x_6, x_3x_6, x_5x_6\}.$$
\noindent
It admits the prime decomposition
$$I=(x_1x_2+x_3x_4,\ x_6)\cap(x_1,\ x_3,\ x_5),$$
\noindent
so that $I$ is reduced and height\,$I=2$. We can apply Lemma \ref{lemma2} to
\begin{eqnarray*}
P_0&=&\{x_1x_6\}\\
P_1&=&\{x_3x_6\}\\
P_2&=&\{x_1x_2+x_3x_4,\ x_5x_6\},
\end{eqnarray*}
\noindent
since $(x_1x_2+x_3x_4)x_5x_6\in (x_1x_6,\ x_3x_6)$. Thus
$$I=\sqrt{(x_1x_6,\ x_3x_6,\ x_1x_2+x_3x_4+x_5x_6)},$$
\noindent
whence ara\,$I\leq 3$. In fact equality holds, so that, in particular, $I$ is not a s.t.c.i.. We show that
\begin{equation}\label{2}3\leq\mbox{ara}\,I\end{equation}
\noindent 
  by means of the well-known inequality (see \cite{Hu}, Theorem 3.4, or \cite{Ha}, Example 2, pp.~414--415) 
$$\mbox{cd}\,I\leq\mbox{ara}\,I.$$
\noindent
Here 
$${\rm cd}\,I=\max\{i\vert H^i_I(R)\ne0\}$$
\noindent 
is called the {\it cohomological dimension} of $I$; $H^i_I$ denotes the $i$-th local cohomology group with respect to $I$.  We refer to Huneke \cite{Hu}  for the basic notions. Our claim (\ref{2}) will follow once we have proven that 
\begin{equation}\label{3} 3\leq\mbox{cd}\,I.\end{equation}
\noindent
According to \cite{Hu}, Theorem 2.2, we have the following long exact sequence of local cohomology groups:
\begin{equation}\label{4} 
\cdots\rightarrow H^3_{I+(x_6)}(R)\rightarrow H^3_I(R)\rightarrow H^3_I(R_{x_6})\rightarrow  H^4_{I+(x_6)}(R)\rightarrow\cdots,\end{equation}
\noindent
where, by \cite{Hu}, Proposition 1.10, 
\begin{equation}\label{5}H^3_I(R_{x_6})\simeq H^3_{I_{x_6}}(R_{x_6}).\end{equation}
\noindent Now
\begin{eqnarray*} I+(x_6)&=&(x_1x_2+x_3x_4,\ x_6)\subset R,\\
I_{x_6}&=&(x_1,\ x_3,\ x_5)\subset R_{x_6},
\end{eqnarray*}
\noindent
and both ideals are generated by a regular sequence, of length 2 and 3 respectively. According to  \cite{Ha}, Example 2, pp.~414--415, it follows that cd\,$(I+(x_6))=2$, and cd\,$I_{x_6}=3$, whence
\begin{eqnarray} H^i_{I+(x_6)}(R)&=&0\qquad\mbox{for }i>2,\label{cd1}\\
H^3_{I_{x_6}}(R_{x_6})&\neq&0.\label{cd2}
\end{eqnarray}
 \noindent
In view of (\ref{5}), (\ref{cd1}) and (\ref{cd2}), from (\ref{4}) we derive an exact row
$$\begin{array}{ccccccc}
0&\rightarrow& H^3_I(R)&\rightarrow& H^3_I(R_{x_6})&\rightarrow&  0\\
&&&&\not{\|}&\\
&&&&0&
\end{array},
$$
\noindent
which implies that 
$$H^3_I(R)\neq0.$$
\noindent
Our claim (\ref{3}) follows.
}\end{example}
Note that Lemma \ref{lemma1}, the original version of Schmitt and Vogel's result, would not allow us to compute ara\,$I$, since it only gives the trivial upper bound
 ara\,$I\leq 4$.\newline
In the next section we shall study an ideal obtained by a slight modification of the ideal $I$ in Example \ref{example1}, namely $J=(x_1x_2+x_3x_4,\ x_1x_5, \ x_3x_5).$
 We shall show that ara\,$J=2$. This result, evidently, cannot be derived from Lemma \ref{lemma2}. A new theorem needs to be introduced. 
\section{A linear algebraic criterion}
Let ${\bf e}_1,\dots, {\bf e}_n$ be the standard basis of the free $R$-module $R^n$. 
\begin{theorem}\label{theorem1}
Let $p_1,\dots, p_{n-1}\in R$ and consider the $n\times (n-1)$-matrix
$$A=\left(\begin{array}{c}
c_1{\bf e}_{i_1}\\
\vdots\\
c_{n}{\bf e}_{i_{n}}
\end{array}
\right),
$$
\noindent
where $c_k\in R$ for $k=1,\dots, n$ and $i_1,\dots, i_{n}\in\{1,\dots, n-1\}$. Let 
\begin{equation}\label{p0}p_0=\sum_{k=1}^{n}(-1)^k\alpha_{k0}\Delta_k,\end{equation}
\noindent
 where $\alpha_{k0}\in R$ for all $k=1,\dots, n$ and $\Delta_k$ is the determinant of the $(n-1)\times (n-1)$ matrix obtained by dropping the $k$-th row of $A$. For all $k=1,\dots, n$ set
\begin{equation}\label{10'}q_k=\alpha_{k0}p_0+c_kp_{i_k}.\end{equation}
\noindent
Then
$$\sqrt{(p_0,c_1p_{i_1},\dots, c_{n}p_{i_{n}})}=\sqrt{(q_1,\dots, q_{n})}.$$
\end{theorem}
\demo It suffices to show that 
$$p_0,c_1p_{i_1},\dots, c_{n}p_{i_{n}}\in\sqrt{(q_1,\dots, q_{n})}.$$
\noindent
For convenience of notation, let us write 
$$A=(\alpha_{ij})_{{\scriptstyle i=1,\dots, n}\atop{\scriptstyle j=1,\dots, n-1}},$$
\noindent
  so that, for $k=1,\dots, n$,
$$q_k=\alpha_{k0}p_0+\sum_{j=1}^{n-1}\alpha_{kj}p_j.$$
\noindent
Then
\begin{eqnarray*} \sum_{k=1}^{n}(-1)^k\Delta_kq_k&=&
\left(\sum_{k=1}^{n}(-1)^k\Delta_k\alpha_{k0}\right)p_0+\sum_{j=1}^{n-1}\left(\sum_{k=1}^{n}(-1)^k\Delta_k\alpha_{kj}\right)p_j\\&=&
p_0^2.
\end{eqnarray*}
\noindent
Here we used (\ref{p0}) and the fact that $\sum_{k=1}^{n}(-1)^k\Delta_k\alpha_{kj}$ is zero: it is the Laplace expansion, with respect to the $j$-th column, of the determinant of the $n\times n$ matrix obtained by adding a copy of the $j$-th column to $A$ on the left. Thus 
$p_0^2\in(q_1,\dots, q_{n}),$
 and $$p_0\in\sqrt{(q_1,\dots, q_{n})}.$$
\noindent In view of (\ref{10'}) we conclude that, for all $k=1,\dots, n$, 
$$c_kp_{i_k}=q_k-\alpha_{k0}p_0\in\sqrt{(q_1,\dots, q_{n})}.$$
\noindent
This completes the proof.
\par\medskip\noindent
The above theorem contains a linear algebraic condition under which an ideal generated by $n+1$ elements is generated, up to radical, by $n$ elements. Analogous criteria were proved in \cite{B2} for the products of two ideals.  Next we derive from Theorem \ref{theorem1} further similar results that apply to ideals which are the sum of a principal ideal and a product ideal. 
\begin{corollary}\label{corollary1} Let $\alpha_1,\alpha_2,\beta_1,\beta_2,\gamma\in R$. Then 
\begin{eqnarray*}&&\sqrt{(\alpha_1\beta_1+\alpha_2\beta_2,\ \beta_1\gamma,\ \beta_2\gamma)}=\\
&&\qquad\qquad\qquad\qquad\qquad\sqrt{(\alpha_1(\alpha_1\beta_1+\alpha_2\beta_2)+\beta_2\gamma,\ \alpha_2(\alpha_1\beta_1+\alpha_2\beta_2)-\beta_1\gamma)}.\end{eqnarray*}
\end{corollary} 
\demo It suffices to apply Theorem \ref{theorem1} for $n=2$, $p_1=\gamma$, $\alpha_{10}=\alpha_1$, $c_1=\beta_2$, $\alpha_{20}=\alpha_2$, $c_2=-\beta_1$. In this case 
$$A=\left(\begin{array}{c}\beta_2\\-\beta_1\end{array}\right),$$
\noindent so that $\Delta_1=-\beta_1$, $\Delta_2=\beta_2$, and
$$p_0=-\alpha_{10}\Delta_1+\alpha_{20}\Delta_2=\alpha_1\beta_1+\alpha_2\beta_2.$$
\par\medskip\noindent
\begin{example}\label{example2}{\rm Let $K$ be a field and consider the ideal
$$J=(x_1x_2+x_3x_4,\ x_1x_5,\ x_3x_5)\in K[x_1,x_2,x_3,x_4,x_5].$$
\noindent
It has a prime decomposition 
$$J=(x_1x_2+x_3x_4,\ x_5)\cap(x_1,\ x_3),$$
\noindent
hence it is reduced of height 2. According to Corollary \ref{corollary1} we have $I=\sqrt{(q_1,q_2)}$ with
\begin{eqnarray}\label{6}
q_1&=&x_2(x_1x_2+x_3x_4)+x_3x_5\\
q_2&=&x_4(x_1x_2+x_3x_4)-x_1x_5.
\end{eqnarray}
\noindent
Here $\alpha_1=x_2$, $\alpha_2=x_4$, $\beta_1=x_1$, $\beta_2=x_3$, $\gamma=x_5$. In particular, $I$ is a s.t.c.i.. 
}
\end{example}
The next result is a generalization of Corollary \ref{corollary1}.
\begin{proposition}\label{proposition1}
Let $\alpha_1, \alpha_2, \beta_1, \beta_2,\gamma_1,\dots, \gamma_{n-1}\in R$ and consider the ideal
$$J=(\alpha_1\beta_1+\alpha_2\beta_2)+(\beta_1, \beta_2)(\gamma_1,\dots, \gamma_{n-1})\subset R.$$
\noindent
Then there are $q_1,\dots, q_n\in (\beta_1, \beta_2)$ such that $$\sqrt{J}=\sqrt{(q_1,\dots, q_n)}.$$
\end{proposition}
\demo We proceed by induction on $n\geq2$. For $n=2$ the claim is Corollary \ref{corollary1}. Let $n>2$ and suppose that the claim is true for $n-1$. Then there are $q_1,\dots, q_{n-1}\in (\beta_1, \beta_2)$, such that 
\begin{equation}\label{induction}\sqrt{(\alpha_1\beta_1+\alpha_2\beta_2)+(\beta_1, \beta_2)(\gamma_1,\dots, \gamma_{n-2})}=\sqrt{(q_1,\dots, q_{n-1})}.\end{equation}
\noindent
Let $\alpha'_1, \alpha'_2\in R$  be such that
\begin{equation}\label{q1} q_1=\alpha'_1\beta_1+\alpha'_2\beta_2.\end{equation}
\noindent  
Then set
\begin{eqnarray}
q'_1&=&\alpha'_1q_1+\beta_2\gamma_{n-1}\label{q'1}\\
q_n&=&\alpha'_2q_1-\beta_1\gamma_{n-1}\label{qn}
\end{eqnarray}
\noindent
Then, $q'_1, q_n\in (\beta_1, \beta_2)$, and, by Corollary \ref{corollary1}, 
\begin{equation}\label{two}\sqrt{(q_1,\beta_1\gamma_{n-1}, \beta_2\gamma_{n-1})}=\sqrt{(q'_1, q_n)},\end{equation}
\noindent so that, in particular, $q_1\in\sqrt{(q'_1,q_n)}$.
 Hence, by (\ref{induction}) and (\ref{two}), 
\begin{eqnarray*}
&&\!\!\!\!\!\!\!\!\sqrt{(q'_1, q_2,\dots, q_n)}=\sqrt{(q'_1,q_1, q_2,\dots, q_n)}\\
&&\qquad\supset\sqrt{(q'_1,q_n)}+\sqrt{(q_1,q_2,\dots, q_{n-1})}\\
&&\qquad=\sqrt{(q_1, \beta_1\gamma_{n-1}, \beta_2\gamma_{n-1})}+
\sqrt{(\alpha_1\beta_1+\alpha_2\beta_2)+(\beta_1, \beta_2)(\gamma_1,\dots, \gamma_{n-2})}\\
&&\qquad\supset(\beta_1\gamma_{n-1}, \beta_2\gamma_{n-1})+(\alpha_1\beta_1+\alpha_2\beta_2)+(\beta_1, \beta_2)(\gamma_1,\dots, \gamma_{n-2})\\
&&\qquad=J.
\end{eqnarray*}
\noindent
On the other hand,  $q_1,\dots, q_{n-1}\in \sqrt J$ by (\ref{induction}) and, consequently, $q_1'\in\sqrt J$ by (\ref{q'1}) and $q_n\in\sqrt J$ by (\ref{qn}). This implies that $\sqrt{(q'_1, q_2,\dots, q_n)}\subset\sqrt J$ and completes the proof.
\par\medskip\noindent
The proof of Proposition \ref{proposition1} yields a recursive construction of $q_1,\dots, q_n$. In Example \ref{example2} we obtained two polynomials $q_1,q_2$ generating the ideal
$$I=(x_1x_2+x_3x_4)+(x_1, x_3)(x_5)$$
\noindent
up to radical. We use this result to construct $q'_1,q_3\in K[x_1,\dots, x_6]$ such that $q_1', q_2, q_3$ generate the ideal
$$J=(x_1x_2+x_3x_4)+(x_1, x_3)(x_5, x_6)$$
\noindent
up to radical. We stick to the notation introduced in Example \ref{example2} and, moreover, we set $\gamma_1=x_5$, $\gamma_2=x_6$. Since, by (\ref{6}),
$$q_1=x_2^2x_1+(x_2x_4+x_5)x_3,$$
\noindent
comparison with (\ref{q1}) yields $\alpha'_1=x_2^2$, $\alpha'_2=x_2x_4+x_5$. Then, according to (\ref{q'1}) and (\ref{qn}), we set
\begin{eqnarray*}
q'_1&=&x_2^2q_1+x_3x_6\\
q_3&=&(x_2x_4+x_5)q_1-x_1x_6.
\end{eqnarray*}
\noindent
The conclusion is: the ideal 
$$J=(x_1x_2+x_3x_4,\ x_1x_5, \ x_3x_5, \ x_1x_6, \ x_3x_6)$$
\noindent
is generated up to radical by 
\begin{eqnarray*}
q'_1&=&x_1x_2^4+x_2^3x_3x_4+x_2^2x_3x_5+x_3x_6\\
q_2&=&x_1x_2x_4+x_3x_4^2-x_1x_5\\
q_3&=&x_1x_2^3x_4+x_2^2x_3x_4^2+x_1x_2^2x_5+2x_2x_3x_4x_5+x_3x_5^2-x_1x_6.
\end{eqnarray*}
\noindent
In general, for all $n\geq5$, the ideal
\begin{eqnarray*}I_n&=&(x_1x_2+x_3x_4,\ x_1x_5, \dots, x_1x_n, \ x_3x_5, \dots, \ x_3x_n)\\
&=& (x_1x_2+x_3x_4,\ x_5, \dots, x_n)\cap (x_1,x_3)\subset K[x_1,\dots, x_n],
\end{eqnarray*}
\noindent
is reduced of height 2 and there are $q_1,q_2, \dots, q_{n-3}\in K[x_1,\dots, x_n]$ generating $I_n$ up to radical. Hence ara\,$I_n\leq n-3$. We are going to show that equality holds, i.e., that the variety $V=V(I_n)\subset K^n$ cannot be defined by $n-4$ equations. This time we use the following cohomological criterion, due to Newstead.
\begin{lemma}\label{Newstead}{\rm [\cite{BS}, Lemma 3]} Let
$W\subset\tilde W$ be affine varieties. Let $d=\dim\tilde
W\setminus W$. If there are $s$ equations $F_1,\dots, F_s$ such
that $W=\tilde W\cap V(F_1,\dots,F_s)$, then 
$$\het^{d+i}(\tilde W\setminus W,{\bz}/r{\bz})=0\quad\mbox{ for all
}i\geq s$$ and for all $r\in{\bz}$ which are prime to {\rm char}\,$K$.
\end{lemma}
We refer to \cite{Mi} for the basic notions on \'etale cohomology. 
Let $p$ be a prime such that $p\ne$\,char\,$K$. In view of Lemma \ref{Newstead}, for our purpose it suffices to show that 
\begin{equation}\label{7} \het^{2n-4}(K^n\setminus V, {\bz}/p{\bz})\ne0.
\end{equation}
\noindent
 By Poincar\'e Duality (see \cite{Mi}, Corollary 11.2, p.~276) we have
\begin{equation}\label{8} \mbox{Hom}(\het^{2n-4}(K^n\setminus V, {\bz}/p{\bz}),{\bz}/p{\bz}) \simeq  \hc^{4}(K^n\setminus V, {\bz}/p{\bz}),
\end{equation}
\noindent
where $\hc$ denotes \'etale cohomology with compact support. 
Hence (\ref{7})  is proven once we have shown that
\begin{equation}\label{9}\hc^{4}(K^n\setminus V, {\bz}/p{\bz})\ne0
\end{equation}
\noindent
For the sake of simplicity, we shall omit the coefficient group ${\bz}/p{\bz}$ henceforth. We have a long exact sequence of cohomology with compact support:
\begin{equation}\label{10}\cdots\rightarrow\hc^3(K^n)\rightarrow \hc^3(V)\rightarrow \hc^4(K^n\setminus V)\rightarrow 
\hc^4(K^n)\rightarrow\cdots\end{equation}
\noindent
It is well known that 
\begin{equation}\label{11}\hc^i(K^t)\simeq\left\{\begin{array}{cl} {\bz}/p{\bz}&\mbox{if }i=2t\\
0&\mbox{else, }
\end{array}\right.
\end{equation}
\noindent
and 
\begin{equation}\label{12}\hc^i(K^t\setminus\{0\})\simeq\left\{\begin{array}{cl} {\bz}/p{\bz}&\mbox{if }i=1,2t\\
0&\mbox{else. }
\end{array}\right.
\end{equation}
\noindent
In particular, being $n\geq5$, we have $\hc^3(K^n)=\hc^4(K^n)=0$, i.e.,  in (\ref{10}) the leftmost and the rightmost groups are zero. It follows that the two middle groups are isomorphic:
$$\hc^3(V)\simeq\hc^4(K^n\setminus V).$$
\noindent
Thus our claim (\ref{9}) is equivalent to
\begin{equation}\label{13} \hc^3(V)\neq0\end{equation}
\noindent
Let $W$ be the subvariety of $K^n$ defined by $x_1=x_3=0$.  
\noindent
Then  $W$ is a $(n-2)$-dimensional affine space over $K$, $W\subset V$, and 
\begin{eqnarray*}
&&\!\!\!\!\!\!\!\!V\setminus W=\\
&&\!\!\!\!\!\!\!\!\{(x_1,\dots, x_n)\in K^n\vert x_1x_2+x_3x_4=x_5=x_6=\cdots =x_n=0\}\setminus\\
&&\qquad\qquad\qquad\qquad\qquad\qquad\qquad\qquad\qquad\{(x_1,\dots, x_n)\in K^n\vert x_1=x_3=0\}.
\end{eqnarray*}
\noindent
Hence $V\setminus W$ can be viewed as a hypersurface $H$ in the 4-dimensional affine space over $K$ (in the coordinates $x_1, x_2, x_3, x_4$) defined by $x_1x_2+x_3x_4=0$, minus the 2-dimensional linear subspace $L$ defined by $x_1=x_3=0$. We have the following diagram with exact rows
\begin{equation}\label{14}
\begin{array}{ccccccc}
\hc^2(W)&\rightarrow&\hc^3(V\setminus W)&\rightarrow&\hc^3(V)&\rightarrow&\hc^3(W)\\
&&\|&&&\\
\hc^2(L\setminus\{0\})&\rightarrow&\hc^3(H\setminus L)&\rightarrow&\hc^3(H\setminus\{0\})&\rightarrow&\hc^3(L\setminus\{0\})\end{array}
\end{equation}
\noindent
Since $W\simeq K^{n-2}$, and $n-2\geq 3$, from (\ref{11}) it follows that $\hc^2(W)=\hc^3(W)=0$; from (\ref{12}) it follows that $\hc^2(L\setminus\{0\})=\hc^3(L\setminus\{0\})=0$. Hence the two middle maps in (\ref{14}) are isomorphisms, so that
\begin{equation}\label{15}\hc^3(V)\simeq \hc^3(H\setminus\{0\}) \end{equation}
\noindent
Note that, up to changing the sign of $x_3$, $\bar H=H\setminus\{0\}$ is the set of all  non zero $2\times 2$ matrices $(x_1, x_3; x_4, x_2)$ with entries in $K$ having proportional rows. The set of such matrices where the first row is zero is a closed subset of $\bar H$ which can be identified with $Z=K^2\setminus\{0\}$, and its complementary set is $\bar H\setminus Z\simeq K^2\setminus\{0\}\times K$.  We thus have a
long exact sequence of \'etale cohomology with compact support:
\begin{equation}\label{sequence}\cdots\rightarrow\hc^2(Z)\rightarrow \hc^3(\bar H\setminus Z)\rightarrow \hc^3(\bar H)\rightarrow \hc^3(Z)\rightarrow\cdots,\end{equation}
\noindent
where, by (\ref{12}), $\hc^2(Z)=\hc^3(Z)=0$, and, moreover, by the K\"unneth formula (\cite{Mi}, Theorem 8.5, p.~258), 
$$\hc^3(\bar H\setminus Z)\simeq\bigoplus_{i+j=3}\hc^i(K^2\setminus\{0\})\otimes\hc^j(K)\simeq \hc^1(K^2\setminus\{0\})\otimes\hc^2(K)\simeq \bz/p\bz,$$
where we again used (\ref{11}) and (\ref{12}).
\noindent
It follows that (\ref{sequence}) gives rise to an isomorphism
\begin{equation}\label{16}\hc^3(\bar H)\simeq{\bz}/p{\bz}.\end{equation}
\noindent
Now (\ref{16}) and (\ref{15}) imply (\ref{13}). This proves our claim (\ref{9}). We have thus shown that 
$$\mbox{ara}\,I_n=n-3.$$
\noindent
It follows that $I_n$ is not a s.t.c.i. for $n>5$.


\begin{thebibliography}{BMT}
\bibitem{B0} Barile, M., Arithmetical ranks of ideals associated to symmetric and alternating matrices. {\it J. Algebra}, {\bf 176}, (1995), 59--82. 
\bibitem{B1} Barile, M., On the number of equations defining certain varieties. {\it Manuscripta Math.}, {\bf 91}, (1996), 483--494.
\bibitem{B2} Barile, M., On the computation of arithmetical ranks. {\it Int.~J.~Pure Appl.~Math.}, {\bf 17}, (2004), 143--161.
\bibitem{B3} Barile, M., On ideals whose radical is a monomial ideal. Preprint (2003). To appear in: {\it Comm.~Algebra}.
\bibitem{B4} Barile, M., Certain minimal varieties are set-theoretic complete intersections. Preprint (2005), arxiv:math.AG/0509475. 
\bibitem{BMT} Barile, M., Morales, M., Thoma, A., Set-theoretic complete intersections on binomials. {\it Proc.~Amer.~Math.~Soc.},  {\bf 130}, (2002), 1893--1903.
\bibitem{Br} Bresinsky, H., Monomial Gorenstein curves in ${\bf A}^4$ are set-theoretic complete intersections, {\it Manuscripta Math.}, \textbf{27}, (1979), 353--358.
\bibitem{BS} Bruns, W.; Schw\"anzl, R., The number of
equations defining a determinantal variety, {\it Bull.~London Math.~Soc.}, {\bf 22}, (1990), 439--445.
\bibitem{El} Eliahou, S. Id\'eaux de d\'efinition des courbes monomiales, in: {\it Complete Intersections}, (Greco, S.; Strano, R. eds.), Lect.~Notes Math.~\textbf{1092} Springer, Heidelberg, 1984, 229--240.
\bibitem{E} Eto, K., Almost complete intersection monomial curves in ${\bf A}^4$, {\it Comm.~Algebra}, \textbf{22}, (1994), 5325--5342.
\bibitem{Ha} Hartshorne, R., Cohomological dimension of algebraic varieties, {\it Ann.~of Math.}, \textbf{88}, (1989), 403--450.
\bibitem{Hu} Huneke, C., Lectures on local cohomology (with an appendix by Amelia Taylor). Available at \underline{http://www.math.ku.edu/\~{}huneke/Vita/chicago-lc.pdf}. 
\bibitem{Mi} Milne, J., {\it \'Etale Cohomology}. Princeton University Press, Princeton, 1980.
\bibitem{RV} Robbiano, L., Valla, G., Some curves in ${\bf P}^3$ are set-theoretic complete intersections, in: {\it Algebraic Geometry-Open problems},  (Ciliberto, C.; Ghione, F.; Orecchia, F., eds.),
Lect.~Notes Math.~{\bf 997}, Springer, Berlin-Heidelberg-New York-Tokyo, 1983, 391--399. 
\bibitem{SV} Schmitt, Th.; Vogel, W., Note on Set-Theoretic Intersections of Subvarieties of Projective Space, {\it Math.~Ann.}, {\bf 245}, (1979), 247--253.
\end{thebibliography}
\end{document}